\documentstyle[12pt]{article}
\input amssym.def
\input amssym.tex 

\newenvironment{demo}[1]%
{\vskip-\lastskip\medskip
  \noindent
  {\em #1.}\enspace
  }%
{\qed\par\medskip
  }

\newcommand{\qed}{
  \strut\hfill
  \mbox{$\Box$}
  }

\newtheorem{theorem}{Theorem}

\newtheorem{lemma}{Lemma}

\newtheorem{corollary}{Corollary}

\newcommand{\hf}{
  \frac12
  }

\begin{document}
\title{Dimension of a minimal nilpotent orbit }
\author{
  Weiqiang Wang\thanks{
  Partially supported by NSF grant DMS-9304580.
       }
}
\date{}
\maketitle
\begin{abstract}{We show that the dimension of the 
minimal nilpotent  coadjoint orbit
for a complex simple Lie algebra 
is equal to twice the dual Coxeter number minus two.}
\end{abstract}

Let $\frak g$ be a finite dimensional complex simple Lie algebra.  
We fix a Cartan subalgebra $\frak h$, a root 
system $\Delta \subset {\frak {h}} ^{*} $ and a set of positive roots
$\Delta _{+} \subset \Delta$.
Let $\rho$ be half the sum of all positive roots.
Denote by $\theta$ the 
highest root and normalize the Killing form 
$$ (\,,\,): \frak g \times \frak g \rightarrow \Bbb C $$
by the condition $(\theta ,\theta ) = 2 $. The dual Coxeter number 
$h^{\vee}$ can be defined as $h^{\vee} = (\rho, \theta) +1 $ (cf. \cite{K}). 
This intrinsic number of the Lie algebra $\frak g$
plays an important role in representation theory(cf. e.g. \cite{K}).

As is well known there exists a unique nonzero
nilpotent (co)adjoint orbit of 
minimal dimension. A coadjoint orbit can be identified with an
adjoint one by means of the Killing form. For more detail on
the nilpotent orbits, we refer to the excellent
exposition \cite{CM} and the references therein. Our result of
this short note is the following theorem.

\begin{theorem}
  The dimension of the minimal nonzero nilpotent orbit equals
  $2h^{\vee} - 2$.
\end{theorem}

We start with the following well-known lemma, cf. for example,
Lemma 4.3.5, \cite{CM}.

\begin{lemma} 
  The dimension of the minimal nonzero nilpotent orbit is equal to
  one plus the number of positive roots not orthogonal to $\theta$.
  
  \label{lemma_1}
\end{lemma}

We call a root $\alpha$ in $\Delta_{+}$ {\em special} if $\theta - \alpha$ 
is also a root. The subset of special roots, denoted by 
$\Bbb S\/$, was singled out in \cite{KW, W} for some 
other purposes. It is easy to see that
we can also define the set $\Bbb S\/$ equivalently
as follows.

\begin{lemma}
  The set $\Bbb S$ is characterized by the property:
  $r_{\theta} (\alpha) = \alpha - \theta, $ if $\alpha \in \Bbb
  S$; $r_{\theta} (\alpha) = \alpha $, if $\alpha \in
  \Delta_{+}-(\Bbb S \cup \{\theta \})$. In other words,
  $\Bbb S \cup \{\theta \}$ is the set of positive roots not orthogonal 
  to $\theta$.

  \label{lemma_2}
\end{lemma}

The following lemma is taken from \cite{KW, W}. The simple proof 
given below follows \cite{W}.
\begin{lemma}   
  The number of special roots is $\# \Bbb S = 2(h^{\vee}-2)$.
  \label{lemma_3}
\end{lemma}

\begin{demo}{Proof}
Since $(\theta, \theta) =2$ and $ (\rho, \theta) = h^{\vee} -1,$ we have
 \begin{equation}
 \label{eqn:num1}
   r_{\theta} \rho = \rho - \frac{2(\rho, \theta)}{(\theta, \theta)} \theta
= \rho - (h^{\vee} -1) \theta.
 \end{equation}

On the other hand, it follows from Lemma \ref{lemma_2} that
 \begin{eqnarray}
    r_{\theta} \rho &= & r_{\theta} 
   \left(\frac{1}{2}\sum_{\alpha \in \Delta_{+}}
   \alpha \right) \nonumber \\
   & = & \hf \left(
                \sum_{{\theta} \neq\alpha \in \Delta_{+}} r_{\theta} (\alpha)
                  - \theta 
             \right) \nonumber \\
   & = & \hf \left(
                \sum_{{\theta} \neq\alpha \in \Delta_{+}} \alpha
                - (\# \Bbb S) \theta - \theta  
             \right) \nonumber \\
   & = & \rho - \hf \left(
                       \# \Bbb S +2 
                    \right) \theta.
\label{eqn:num2}
\end{eqnarray}
Thus this lemma follows by comparing the right hand sides of the equations
(\ref{eqn:num1}) and (\ref{eqn:num2}).
\end{demo}

By combining Lemmas \ref{lemma_1}, \ref{lemma_2} and \ref{lemma_3}, we
prove our theorem. We have an immediate corollary from 
Lemmas \ref{lemma_2} and \ref{lemma_3}.

\begin{corollary}
  The length of the reflection $r_{\theta}$ is 
  $l \left( r_{\theta} \right) = 2h^{\vee} -3\/$.
\end{corollary}

\noindent
Max-Planck Institut f\"ur Mathematik, 53225 Bonn, Germany\\
E-mail address: {\tt wqwang@mpim-bonn.mpg.de}

\end{document}